# The Radical Mistakes in the Theories on Calculus of Cauchy-Lebesgue System

Xiaoping Ding
*Beijing institute of pharmacology of 21st century*

**Abstract** The misunderstanding of the concept of differentials in the theories on calculus of Cauchy-Lebesgue system was exposed in this paper. The defects of the definition of differentials and the associated mistakes in the differentiation of composite functions were pointed out and discussed.

**Introduction**

In this paper it is suggested that the history of calculus consists of three stages including the period of precalculus before I. Newton established his calculus system initially in 1667, the first period of history of calculus from 1667 until A. L. Cauchy published his monograph 'cours d'analyse' in 1821, the second period of history of calculus from 1821 until the novel calculus system was established in 2010 and the third one after 2010.

During the period of 189 years from 1821 to 2010, the theories on calculus developed in this second period has never been justified by the application of mathematics in science and technology, which only demonstrated the usefulness of calculus invented by I. Newton and G. W. Leibniz and developed by the subsequent mathematicians as before. In contrast, the history of 189 years has revealed that the theories on calculus in this second period have blocked the development of mathematics significantly.

In this paper the radical mistakes in the theories on calculus of Cauchy-Lebesgue system were shown and discussed. The preliminary achievements of novel theories on calculus established by the author have been published[1].

**I The misunderstanding of the concept of differentials by Cauchy**

The basic ideas of theories on calculus in the first period (i.e., Newton-Leibniz theories on calculus) are correct although they are not consistent. Cauchy agreed with formulas $dx^e = ex^{e-1}dx$, $\int x^e dx = \dfrac{x^{e+1}}{e+1}$ and $\int_a^b ydx = Z(b) - Z(a)$, given by Leibniz, but he did not understand Leibniz's explanation of the differentials $dx$ and $dy$.

Leibniz said, 'A differential is like the contact angle of Euclid, which is smaller than any given quantity, but not equal to zero.' He also said, 'We consider an infinitesimal quantity as a relative zero, not a simple zero or an absolute zero.' It should be noted that a worldwide mistake exists that the differential is considered as an arbitrarily small quantity. Indeed, Leibniz indicated that the differential is not zero, or a finite quantity, not to mention infinity, but 'a relative zero smaller than any given quantity'.

It is a new type of quantity since the concept of modern numbers hasn't appeared in Leibniz's time; while the concept of 'the relative zero' existed in the concepts of all the numbers, including unnormalized numbers of Robinson, from the time of Cantor-Dedekind to 2010. However, Cauchy introduced the idea of the limit to calculus system because he wasn't able to understand the ideas of Leibniz. Making use of the concept of the limit and the derivative to define the differential violates the original meaning of the differential from Leibniz.

**II The defects of the definition of the differential and the associated mistakes**

Cauchy defined the differential as a finite quantity, which had changed Leibniz's idea that the differential was defined as 'a qualitative zero', 'smaller than any given quantity'. It seems that a derivative has been defined perfectly and the Newton-Leibniz equation has been proved by making use of the concept of limit, but there exist radical problems of equating $dx$ and $\Delta x$ when it comes to defining the differential, because the differential (expressed starting with the symbol 'd') has been defined as the linear main-part of a change (expressed starting with the symbol '$\Delta$') resulting in that $dx$ is not equal to $\Delta x$. Formulas $\Delta z = E'(y)\Delta y + o_E(\Delta y)$, $\Delta y = F'(x)\Delta x + o_F(\Delta x)$ and $\Delta x = G'(t)\Delta t + o_G(\Delta t)$ are obtained regarding arbitrarily differentiable functions $z = E(y)$, $y = F(x)$ and $x = G(t)$. Infinitesimals of higher order $o_E(\Delta y)$, $o_F(\Delta x)$ and $o_G(\Delta t)$ cannot be removed if differentials $dz$, $dy$ and $dx$ are not defined as the linear main-parts of changes $\Delta z$, $\Delta y$ and $\Delta x$; while inequations $dz \neq \Delta z$, $dy \neq \Delta y$ and $dx \neq \Delta x$ are obtained if the differentials are defined as the linear main-parts of the changes, resulting in that $dx$ is not equal to $\Delta x$. Perhaps Cauchy realized that his idea of reestablishing theories on calculus did not make sense, however, he continued to establish his 'theories on calculus' producing more associated mistakes.

The current theories on calculus (i.e., Cauchy-Lebesgue theories on calculus) were 'established' via a nonmathematical way on the definition of the differential based on two assumptions without logic, both presented in the current textbooks on mathematical analysis (calculus). Given a differentiable function $y = F(x_1, x_2, ..., x_n)$, one assumption is that differentials of $x_1, x_2, ..., x_n$ are considered as the differentials of functions $y = x_1, x_2, ..., x_n$, respectively, without logic[2-4]; the other one is that equations $dx_1 = \Delta x_1$, $dx_2 = \Delta x_2$, ..., $dx_n = \Delta x_n$ are defined without logic, where $n$ is equal to an arbitrarily finitely natural number in both cases[4-8].

From a viewpoint of argument, Cauchy has defined the differential (e.g., $dy$ and $dx$) as the linear main-part of a change (e.g., $\Delta y$ and $\Delta x$), therefore, it is not appropriate to define the differential is equal to the change. On the other side, from the differentiable function $y_i = F_i(x_1, x_2, ..., x_j, ..., x_n)$, the general relationship of $y_i$ and $x_j$ has been determined. Thus the assumption of the equations $y_i = x_1 = x_2 = ... = x_j = ... = x_n$, which has transferred the function of n variables into one variable, is not acceptable. Actually the definition of the equation $dx_j = \Delta x_j$ is equivalent to the assumption that the differential of $x_j$ is equal to the differential of the function $y_i = x_j$. One assumption is incorrect, thus the other one is also incorrect. In a word, both of the two assumptions are wrong.

From a viewpoint of counterargument, regarding generally differentiable functions $z_h = E_h(y_1, y_2, ..., y_i, ..., y_m)$, $y_i = F_i(x_1, x_2, ..., x_j, ..., x_n)$ and $x_j = G_j(t_1, t_2, ..., t_k, ..., t_o)$, equations $\Delta z_h = \sum_{i=1}^{m} \frac{\partial E_h}{\partial y_i} \Delta y_i + o_h(\rho_h)$, $\Delta y_i = \sum_{j=1}^{n} \frac{\partial F_i}{\partial x_j} \Delta x_j + o_i(\rho_i)$ and $\Delta x_j = \sum_{k=1}^{o} \frac{\partial G_j}{\partial t_k} \Delta t_k + o_j(\rho_j) = dx_j + o_j(\rho_j)$ are obtained. Therefore, generally speaking, inequations $dz_h \neq \Delta z_h$, $dy_i \neq \Delta y_i$ and $dx_j \neq \Delta x_j$ are obtained. Equations $z_h = y_i (i=1,2,...,m)$, $y_i = x_j (j=1,2,...,n)$ and $x_j = t_k (k=1,2,...,o)$ are not acceptable due to the relationship of equivalence.

The two assumptions are the most representative in the Cauchy-Lebesgue theories on calculus, which violate the principles of science. However, the defenders of Cauchy-Lebesgue system said, 'The subject of mathematics is just a formal system, not a physical one. The definition of the differential by Cauchy doesn't deal with composite functions, therefore, the Cauchy-Lebesgue system is impeccable.' However, in fact it is not true. Many cases in calculus are related to composite functions including the consistence of the expression of the differential[9], the rules of derivation of composite functions, implicit functions, parametric equations and polar equations, the two types of Substitution Rule and the method of substitution of variables in differential equations.

Actually there exist mathematicians supporting Cauchy-Lebesgue theories on calculus, who have also discovered the mistakes. Some of them have done some corrections silently although the corrected system is still wrong[2,9-11]. Some of them choose to avoid dealing with differentials[12].

**Conclusions**

1. Cauchy has misunderstood the meaning of differentials given by I. Newton, which is considered as an appropriate explanation of the differential.

2. The definition of differentials in the theories on calculus of Cauchy-Lebesgue system is not logical and leads to mistakes in differentiation of composite functions.

**References**

1. Ding, X. A novel system of theories on calculus. *Quality Education* **10**, 71-73 (2011).

2. Хинчин, А. Я. *Eight Lectures on Mathematical Analysis* 79-83 (People's Post Press (Translated in Chinese by Wang, H. and Qi, M.), Beijing, 2010).

3. Mathematics teaching and research section of Tongji University (eds). *Higher Mathematics* (Higher Education Press (in Chinese), Beijing, 1999).

4. Фихтенгольц, Г. М. *Principles of Mathematical Analysis* Vol. 1, 170 (People's Education Press (Translated in Chinese by Wu, Q. and Lu, X.), Beijing, 1959).

5. Zhang, Z. *New Lectures on Mathematical Analysis* 1st edn (Peking University Press (in Chinese), Beijing, 1990).

6. Фихтенгольц, Г. М. *A Course on Differential and Integral Calculus* 8th edn, Vol. 1, 175-329 (Higher Education Press (Translated in Chinese by Yang, T. & Ye, Y.), Beijing, 2006).

7. Никольский, С. М. *A Course on Mathematical Analysis* Vol. 1, 141-155 (People's Education Press (Translated in Chinese by Liu, Y. Guo, S. & Gao, S.), Beijing, 2006).

8. Stewart, J. *Calculus* 6th edn (China People's University Press (Translated in Chinese by Zhang, N.), Beijing, 2009).

9. Архипов, Г. И., Садовничий, В. А. & Чубариков, В. Н. *Lectures of Mathematical Analysis* 3rd edn, 81-255 (Higher Education Press (Translated in Chinese by Wang, K.), Beijing, 2006).

10. Зорич, В. А. *Mathematical Analysis* 159-389 (Higher Education Press (Translated in Chinese by Jiang, D., Wang, K., Zhou, M. & Kuang, R.), Beijing, 2006).

11. Карташев, А. П. & Рождественский, Б. Л. *Mathematical Analysis* 117-301 (Inner Mongolia University Press (Translated in Chinese by Cao, Z. and Ni, X.), Huhehaote, 1991).

12. Rudin, W. *Principles of Mathematical Analysis* 3rd edn (China Machine Press Beijing, 2004).

Correspondence and requests for materials should be addressed to Xiaoping Ding
(e-mail: jxddroc@126.com)